# Exponential functionals
# of Lévy processes*

## Jean Bertoin, Marc Yor

*Laboratoire de Probabilités and Institut universitaire de France,*
*Université Pierre et Marie Curie,*
*175 rue du Chevaleret,*
*F-75013 Paris, France*
*e-mail:* `jbe@ccr.jussieu.fr`*;* `deaproba@proba.jussieu.fr`

**Abstract:** This text surveys properties and applications of the exponential functional $\int_0^t \exp(-\xi_s)ds$ of real-valued Lévy processes $\xi = (\xi_t, t \geq 0)$.

## 1. Introduction

Throughout this work, we shall consider a real-valued **Lévy process** $\xi = (\xi_t, t \geq 0)$, refering to [5, 60] for background. This means that the process $\xi$ starts from $\xi_0 = 0$, has right-continuous paths with left-limits, and if $(\mathcal{F}_t)_{t \geq 0}$ denotes the natural filtration generated by $\xi$, then the increment $\xi_{t+s} - \xi_t$ is independent of $\mathcal{F}_t$ and has the same law as $\xi_s$ for every $s, t \geq 0$. It is well-known (see e.g. Sato [60]) that the distribution of $\xi$ is determined by its one-dimensional marginals, and thus by its characteristic function which has the form

$$\mathbb{E}\left(\exp(\mathtt{i}\lambda\xi_t)\right) = \exp(-t\Psi(\lambda)), \qquad \lambda \in \mathbb{R}, t \geq 0.$$

The characteristic exponent $\Psi$ is given by the celebrated **Lévy-Khintchine formula**

$$\Psi(\lambda) = \mathtt{i}\mathtt{a}\lambda + \frac{1}{2}\sigma^2\lambda^2 + \int_{\mathbb{R}} \left(1 - e^{\mathtt{i}\lambda x} + \mathtt{i}\lambda x \mathbf{1}_{\{|x|<1\}}\right)\Pi(dx), \qquad (1)$$

where $\mathtt{a} \in \mathbb{R}$, $\sigma^2 \geq 0$ is known as the Gaussian coefficient, and $\Pi$ is a positive measure on $\mathbb{R}\backslash\{0\}$ such that $\int(1 \wedge |x|^2)\Pi(dx) < \infty$ which is called the **Lévy measure**. In turn, $\mathtt{a}$, $\sigma^2$ and $\Pi$ are uniquely determined by $\Psi$. As an example, for $\xi_t = \sigma B_t + bt$ where $(B_t, t \geq 0)$ is a standard Brownian motion, we have $\Psi(\lambda) = -\mathtt{i}\lambda b + \frac{1}{2}\sigma^2\lambda^2$, so $\mathtt{a} = -b$ and $\Pi = 0$. In the sequel, the trivial case when $\xi \equiv 0$ (i.e. $\Psi \equiv 0$) will be implicitly excluded.

---

*This is an original survey paper.

In this text, we will be interested in the **exponential functional**

$$I_t := \int_0^t \exp(-\xi_s) ds\,, \qquad t \geq 0$$

(note the $-$ sign in the exponential) and especially in its terminal value

$$I_\infty := \int_0^\infty \exp(-\xi_t) dt\,,$$

under conditions on $\xi$ which ensure the finiteness of $I_\infty$. Our motivation stems from the fact that this functional arises in a variety of settings (self-similar Markov processes, mathematical finance, random processes in random environment, Brownian motion on a hyperbolic space, ...), see the forthcoming Sections 5 and 6. The present paper is intended as a survey of this area; it does not contain new results, and only a selection of the statements are carefully proven (of course we shall provide precise references for the remaining ones).

The rest of this survey consists of five sections, which are devoted to the following topics: Section 2 provides necessary and sufficient conditions for the finiteness of $I_\infty$. In Section 3, we give, under certain assumptions on the Lévy process, explicit formulae for the entire moments of the exponential functional, which often determine uniquely its distribution. In the general case, this distribution can be rather complicated; nonetheless, in Section 4, we discuss some instances for which it can be described explicitely, in particular when $\xi$ is a multiple of the Poisson process. In some cases also, the law of $I_\infty$ is absolutely continuous with respect to Lebesgue measure, and its density satisfies an integro-differential equation. Decompositions of $I_\infty$ before and after a stopping time give rise to solutions of random affine equations. Self-similar Markov processes, which form an important family of Markov processes on $]0,\infty[$ having close connections with the exponential functional, are discussed in Section 5. Finally, we briefly present in Section 6 several areas where results on the exponential functionals or self-similar Markov processes can be fruitfully applied.

## 2. Finiteness of the terminal value

As a preliminary to the study of the exponential functional, we first investigate when its terminal value is finite. The following result relates the finiteness of $I_\infty$ to pathwise properties of the Lévy process $\xi$, which in turn are specified in terms of the one-dimensional distributions of $\xi$ or of its characteristics.

**Theorem 1** *The following assertions are equivalent:*
(i) $I_\infty < \infty$ *a.s.*
(ii) $\mathbb{P}(I_\infty < \infty) > 0$.
(iii) $\xi$ *drifts to* $+\infty$, *i.e.* $\lim_{t\to\infty} \xi_t = +\infty$ *a.s.*
(iv) $\lim_{t\to\infty} t^{-1}\xi_t > 0$ *a.s.*
(v) $\int_1^\infty \mathbb{P}(\xi_t \leq 0) t^{-1}\, dt < \infty\,.$

(vi) *Either*

$$\int_{]-\infty,-1[} |x| \Pi(dx) < \infty \quad and \quad -\mathtt{a} + \int_{|x|>1} x\Pi(dx) \in ]0,\infty], \qquad (2)$$

*or*

$$\int_{]-\infty,-1[} |x| \Pi(dx) = \int_{]1,\infty[} x\Pi(dx) = \infty \quad and \quad \int_1^\infty \overline{\Pi}^-(x) d\left(x/J^+(x)\right) < \infty,$$

*where for every* $x > 0$

$$\overline{\Pi}^+(x) = \Pi(]x,\infty[) \quad , \quad \overline{\Pi}^-(x) = \Pi(]-\infty,-x[) \quad , \quad J^+(x) = \int_0^x \overline{\Pi}^+(y) dy.$$

**Remark.** The first part of condition (2) is equivalent to requiring that the negative part of $\xi_1$, $\xi_1^- := |\xi_1| \mathbf{1}_{\{\xi_1 < 0\}}$, has a finite mean. Then $\xi_1 = \xi_1^+ - \xi_1^-$ has a well-defined (possibly infinite) mean

$$m := \mathbb{E}(\xi_1) = \mathbb{E}(\xi_1^+) - \mathbb{E}(\xi_1^-) = -\mathtt{a} + \int_{|x|>1} x\Pi(dx),$$

and the second part of (2) amounts to requiring that $m > 0$.

**Proof:** It is plain that (iv) $\Rightarrow$ (i) $\Rightarrow$ (ii). That (ii) $\Rightarrow$ (i) follows from the 0-1 law for $\xi$. The equivalence (iii) $\Leftrightarrow$ (v) is the criterion of Rogozin (see Theorem VI.12 in [5]), whereas (iii) $\Leftrightarrow$ (iv) $\Leftrightarrow$ (vi) is a version of Erickson's theorem for Lévy processes (see [6]). Finally, it is known (see Theorem VI.12 in [5]) that when (iii) fails, $\xi$ visits the negative half-line at arbitrarily large times a.s. It is easily seen from an application of the strong Markov property that then $I_\infty = \infty$ a.s., which entails (ii) $\Rightarrow$ (iii) and thus completes the proof. ∎

## 3. Calculation of entire moments

In this section, we are interested in the calculation of the entire moments of the exponential functional evaluated at some independent exponential time $T$, viz. $\mathbb{E}(I_T^k)$ for $k \in \mathbb{Z}$. More precisely, positive entire moments are computed when the Lévy process $\xi$ is a subordinator, whereas negative entire moments are being dealt with when the Lévy process has finite positive exponential moments.

### 3.1. Positive moments in the case of subordinators

In this section, we suppose that $\xi = (\xi_t, t \geq 0)$ is a **subordinator**, i.e. $\xi$ has non-decreasing sample paths. Equivalently, this means that the Gaussian coefficient $\sigma^2$ is equal to 0 and that the Lévy measure does not charge $]-\infty,0]$ and fulfills $\int (1 \wedge x)\Pi(dx) < \infty$, and finally

$$\mathtt{d} := -\mathtt{a} - \int_{[0,1]} x\Pi(dx) \geq 0.$$

The coefficient d above is known as the **drift** of the subordinator.

In this setting, it is more convenient to use Laplace transforms rather than characteristic functions to specify distributions. More precisely, the law of the subordinator $\xi$ is characterized via the Laplace transform

$$\mathbb{E}\left(\exp(-q\xi_t)\right) = \exp(-t\Phi(q)), \qquad t, q \geq 0, \tag{3}$$

where $\Phi : [0, \infty[ \to [0, \infty[$ denotes the so-called Laplace exponent (sometimes also referred to as the Bernstein function) of $\xi$. In turn, $\Phi$ is given by an analytic continuation of the characteristic exponent $\Psi$ to the upper half-plane, namely

$$\Phi(q) = \Psi(\mathtt{i}q), \qquad q \geq 0.$$

The Lévy-Khintchine formula now reads

$$\Phi(q) = \mathtt{d}q + \int_{]0,\infty[} (1 - \mathrm{e}^{-qx})\Pi(dx), \qquad q \geq 0.$$

We refer to [7] for background on subordinators.

The first part of the following theorem has been established by Carmona *et al.* [20], the second is an immediate extension of [15].

**Theorem 2** *Suppose that $\xi$ is a subordinator with Laplace exponent $\Phi$. Let $T$ denote a random time which has an exponential distribution with parameter $q \geq 0$ (for $q = 0$, we have $T \equiv \infty$) and is independent of $\xi$.*
(i) *The positive entire moments of the exponential functional evaluated at time $T$ are given in terms of $\Phi$ by the identity*

$$\mathbb{E}\left(I_T^k\right) = \frac{k!}{(q + \Phi(1)) \cdots (q + \Phi(k))}, \qquad \text{for } k = 1, 2, \ldots \tag{4}$$

*As a consequence, $\mathbb{E}(\exp(aI_T)) < \infty$ whenever $a < \Phi(\infty)$, and the distribution of $I_T$ is determined by (4).*
(ii) *For each $q \geq 0$, there exists a unique probability measure $\rho_q$ on $[0, \infty[$ which is determined by its entire moments*

$$\int_{[0,\infty[} x^k \rho_q(dx) = (q + \Phi(1)) \cdots (q + \Phi(k)) \qquad \text{for } k = 1, 2, \ldots$$

*In particular, if $R_q$ is a random variable with law $\rho_q$ that is independent of $\xi$, then we have the identity in distribution*

$$I_T R_q \overset{\mathcal{L}}{=} \mathbf{e}, \tag{5}$$

*where $\mathbf{e}$ denotes a standard exponential variable.*

Of course, one can rephrase formula (4) in terms of Laplace transform for the moments of the exponential functional evaluated at a time $t$, specifically

$$\int_0^\infty \mathbb{E}\left(I_t^k\right) \mathrm{e}^{-qt} dt = \frac{k!}{q(q + \Phi(1)) \cdots (q + \Phi(k))}, \qquad \text{for } k = 1, 2, \ldots$$

which provides a formal characterization of the distribution of $I_t$. However in practice, an explicit inversion of the Laplace transform does not seem to be an easy matter.

**Proof:** (i) We first establish the formula (4) in the case $q = 0$, that is for $T = \infty$ a.s. Set

$$J_t = I_\infty - I_t = \int_t^\infty \exp(-\xi_s)ds$$

for every $t \geq 0$. On the one hand, for every positive real number $r > 0$, we have the identity

$$J_0^r - J_t^r = r \int_0^t \exp(-\xi_s)J_s^{r-1}ds \,. \tag{6}$$

On the other hand, we may express $J_s$ in the form $J_s = \exp(-\xi_s)I'_\infty$, where

$$I'_\infty = \int_0^\infty \exp(-\xi'_u)du \quad \text{and} \quad \xi'_u = \xi_{s+u} - \xi_s \,. \tag{7}$$

From the independence and stationarity of the increments of the Lévy process, we see that $I'_\infty$ has the same law as $I_\infty$ and is independent of $\xi_s$. Plugging this in (6) and taking expectations, we get using (3) that

$$\begin{aligned}
\mathbb{E}\left(I_\infty^r\right)\left(1 - \exp(-t\Phi(r))\right) &= r \int_0^t \exp(-s\Phi(r))\mathbb{E}(I_\infty^{r-1})\,ds \\
&= \frac{r}{\Phi(r)}\left(1 - \mathrm{e}^{-t\Phi(r)}\right)\mathbb{E}\left(I_\infty^{r-1}\right)\,.
\end{aligned}$$

Finally

$$\mathbb{E}\left(I_\infty^r\right) = \frac{r}{\Phi(r)}\mathbb{E}\left(I_\infty^{r-1}\right)\,,$$

and since $\mathbb{E}(I_\infty^0) = 1$, we get the formula of the statement for $q = 0$ by iteration, taking $r = k \in \mathbb{N}$.

Now we treat the general case when the parameter $q$ of the exponential time is finite. Pick $a > 0$ arbitrarily large, and let $N = (N_t, t \geq 0)$ be an independent Poisson process with parameter $q$. So if $T$ denotes the first jump time of the Poisson process $N$, then $T$ has an exponential distribution with parameter $q$ and is independent of $\xi$. Write $\xi_t^{(a)} = \xi_t + aN_t$, so $\xi^{(a)}$ is again a subordinator with Laplace exponent $\Phi^{(a)}(p) = \Phi(p) + q(1 - \mathrm{e}^{-ap})$, $p \geq 0$. Denote also by $I^{(a)}$ the exponential functional of $\xi^{(a)}$. Plainly, $I_\infty^{(a)}$ decreases to $I_T$ as $a$ increases to $\infty$, and since $\lim_{a\to\infty}\Phi^{(a)}(p) = \Phi(p) + q$ for every $p > 0$, we deduce by dominated convergence that

$$\begin{aligned}
\mathbb{E}(I_T^k) &= \lim_{a\to\infty}\mathbb{E}((I_T^{(a)})^k) \\
&= \lim_{a\to\infty}\frac{k!}{\Phi^{(a)}(1)\cdots\Phi^{(a)}(k)} \\
&= \frac{k!}{(q + \Phi(1))\cdots(q + \Phi(k))}\,.
\end{aligned}$$

The rest of the proof is immediate.

(ii) We refer to [15] for the proof of the existence of $\rho$; see also Berg and Duran [4] for an analytical proof. The second part of the statement now follows from (i), since $k!$ is the $k$-th moment of the standard exponential law and the latter is determined by its entire moments. ∎

In Section 4.1, we shall present a couple of examples in which the moment problem (4) can be solved explicitly, yielding expressions for the law of the terminal value of exponential functional $I$. Theorem 2(ii) also invites to look for factorizations of an exponential variable in two independent factors, one of which may be obtained (in distribution) as the exponential functional $I_\infty$. We refer to [13] for several classical such factorizations which fit in this setting.

We conclude this section by mentioning that moments calculations as in Theorem 2(i) can be extended to the situation when $\xi$ is an increasing additive process, i.e. its increments are independent but not necessarily homogeneous. This situation arises naturally in particular in the representation of self-decomposable laws on $\mathbb{R}_+$ (see Sato [59] and Jeanblanc, Pitman and Yor [41]). It also occurs in in non-parametric Bayesian statistics, where the mean of a random distribution chosen from a neutral-to-the-right prior can be represented as the exponential functional of an increasing additive process. See [28] and the references therein for much more on this topic.

### 3.2. Negative moments

Here, we shall suppose that the Lévy process $\xi$ possesses exponential moments of all positive orders, i.e. $\mathbb{E}(\exp(q\xi_t)) < \infty$ for every $t, q \geq 0$. Recall (e.g. from Theorem 25.3 in Sato [60]) that this is equivalent to assuming that the Lévy measure fulfills

$$\int_{[1,\infty[} \mathrm{e}^{qx} \Pi(dx) < \infty \qquad \text{for every } q > 0\,;$$

this condition holds for instance when the jumps of $\xi$ are bounded from above by some fixed number, and in particular when $\xi$ has only negative jumps. In this situation, the characteristic exponent $\Psi$ has an analytic extension to the half-plane with negative imaginary part, and one has

$$\mathbb{E}(\mathrm{e}^{q\xi_t}) = \exp(t\psi(q)) < \infty\,, \qquad t, q \geq 0\,, \tag{8}$$

where

$$\psi(q) := -\Psi(-\mathtt{i}q)\,, \qquad q \geq 0\,.$$

In this direction, note that the assumption that $\xi$ drifts to $+\infty$, is equivalent to

$$m = \mathbb{E}(\xi_1) = \psi'(0+) \in ]0, \infty[\,, \tag{9}$$

see the remark after Theorem 1. The following result appears in [15].

**Theorem 3** *Assume that (8) and (9) hold. Then, for every integer $k \geq 1$ and for every $0 < t \leq \infty$, $\mathbb{E}(I_t^{-k}) < \infty$; furthermore we have*

$$\mathbb{E}\left(I_\infty^{-k}\right) \, = \, m\,\frac{\psi(1)\cdots\psi(k-1)}{(k-1)!}\,,$$

*with the convention that the right-hand side equals $m$ for $k = 1$.*

*If moreover $\xi$ has no positive jumps, then $1/I_\infty$ admits some exponential moments. Hence, the distribution of $I_\infty$ is determined by its negative entire moments.*

**Proof:** To start with, let us first check the finiteness of the negative moments of $I_t$ for $t > 0$, and *a fortiori* for $I_\infty$. In this direction, we observe from (8) that for every $a, q > 0$,

$$\mathbb{P}(\xi_1 \geq a) \, \leq \, \exp(\psi(q) - qa)\,.$$

Next, we have for all sufficiently large $a$ that

$$\mathbb{P}(\inf_{0 \leq t \leq 1} \xi_t < -a) \, \leq \, 1/2\,,$$

and it follows from an application of the strong Markov property at the first passage time of $\xi$ above the level $2a$ that

$$\mathbb{P}(\sup_{0 \leq t \leq 1} \xi_t \geq 2a) \, \leq \, 2\mathbb{P}(\xi_1 \geq a) \, \leq \, 2\exp(\psi(q) - qa)\,.$$

This entails the bound

$$\mathbb{P}(I_1 < \varepsilon^2) \, \leq \, \mathbb{P}(\sup_{0 \leq t \leq 1} \xi_t > 2\log 1/\varepsilon) \, \leq \, 2\exp(\psi(q))\varepsilon^q$$

for small enough $\varepsilon > 0$, and thus ensures that $\mathbb{E}(I_1^{-k}) < \infty$ for every $k > 0$.

Now the proof is similar to that of Theorem 2(i). More precisely, recall the notation

$$J_t \, = \, \int_t^\infty \exp(-\xi_s)ds$$

for every $t \geq 0$ and observe that for every real number $r > 0$, there is the identity

$$J_t^{-r} - J_0^{-r} \, = \, r\int_0^t \exp(-\xi_s)J_s^{-(r+1)}ds\,. \tag{10}$$

As in (7), we write for every fixed $s$

$$J_s \, = \, \exp(-\xi_s)I_\infty'\,,$$

where $I_\infty'$ has the same law as $I_\infty$ and is independent of $\xi_s$. Plugging this in (10) and taking expectations, we get from (8) that

$$\mathbb{E}\left(I_\infty^{-r}\right)(\exp(t\psi(r)) - 1) \, = \, r\int_0^t \exp(s\psi(r))\mathbb{E}(I_\infty^{-(r+1)})\,ds\,,$$

and finally

$$\mathbb{E}\left(I_\infty^{-(r+1)}\right) \;=\; \frac{\psi(r)}{r}\mathbb{E}\left(I_\infty^{-r}\right).\tag{11}$$

Because we know that the negative moments of $I_\infty$ are finite, we may let $r \to 0+$ and apply the theorem of dominated convergence. We get

$$\mathbb{E}(1/I_\infty) \;=\; \lim_{r\to 0+}\frac{\psi(r)}{r} \;=\; \psi'(0+) \;=\; m\,.$$

The formula of Theorem 3 now follows by induction.

Finally, let us check that the distribution of $1/I_\infty$ is determined by its (positive) entire moments when $\xi$ has no positive jumps. In that case, the Lévy measure $\Pi$ does not charge $[0,\infty[$, and the Lévy-Khintchine formula (1) readily yields the bound $\psi(q) = O(q^2)$ as $q \to \infty$. As a consequence, there is some finite constant $c > 0$ such that

$$\mathbb{E}(I_\infty^{-k}) \;=\; m\,\frac{\psi(1)\cdots\psi(k-1)}{(k-1)!} \;\leq\; c^k k! \qquad \text{for every } k = 1,\dots.$$

This entails that

$$\mathbb{E}\left[\exp\left(\frac{1}{2cI_\infty}\right)\right] \;<\; \infty\,,$$

hence the distribution of $1/I_\infty$ is characterized by its entire moments. ∎

We make the conjecture that the positive entire moments of $1/I_\infty$ do not determine its distribution when the Lévy process $\xi$ has positive jumps (it has been shown in [15] that this is indeed the case - discussed in the next subsection - when $\xi$ is a Poisson process, and more precisely, that the negative entire moments of the exponential functional of a Poisson process are closely related to those of the log-normal distribution). We refer to [62, 64] and the references therein for more about indeterminacy for the moment problem.

## 4. On the distribution of the exponential functional

### *4.1. Some important special cases*

In this section, we briefly present some explicit results on the distribution of the exponential functional in three important special cases; further examples can be found e.g. in [20, 21].

• **Standard Poisson process:** Here, we consider a standard Poisson process $(N_t, t \geq 0)$ and define for $q \in ]0,1[$ its exponential functional by

$$I_\infty^{(q)} \;=\; \int_0^\infty q^{N_t}\,dt\,.$$

In other words, $I_\infty^{(q)} = I_\infty$ for $\xi = -(\log q)N$. Note that we may also express $I_\infty^{(q)}$ in the form

$$I_\infty^{(q)} = \sum_{n=0}^{\infty} q^n \varepsilon_n, \tag{12}$$

where $\varepsilon_n = T_{n+1} - T_n$, $n = 0, 1, \ldots$ denote the waiting times between the successive jump times $T_n = \inf\{t : N_t = n\}$ of $(N_t, t \geq 0)$. In other words, $(\varepsilon_n, n \in \mathbb{N})$ is a sequence of i.i.d. exponential variables with parameter 1.

In order to specify the distribution of the exponential functional, it is quite convenient to use the so-called $q$-calculus (see, e.g. [31], [44], ...) which is associated with the basic hypergeometric series of Euler, Gauss, ... More precisely, let us introduce some standard notation in this setting:

$$(a; q)_n = \prod_{j=0}^{n-1} (1 - aq^j) \quad , \quad (a; q)_\infty = \prod_{j=0}^{\infty} (1 - aq^j),$$

and

$$\Gamma_q(x) = \frac{(q; q)_\infty}{(q^x; q)_\infty} (1 - q)^{1-x}.$$

We may now state the main formulas for the distribution of the exponential functional of the Poisson process, which are excerpts from [10] (see also [25]): The Laplace transform of $I_\infty^{(q)}$ is given by

$$\mathbb{E}\left(\exp(\lambda I_\infty^{(q)})\right) = \frac{1}{(\lambda; q)_\infty} \qquad (\lambda < 1), \tag{13}$$

its Mellin transform by

$$\mathbb{E}\left(\left(I_\infty^{(q)}\right)^s\right) = \frac{\Gamma(1+s)}{\Gamma_q(1+s)(1-q)^s} = \Gamma(1+s)\frac{(q^{1+s}; q)_\infty}{(q; q)_\infty}, \tag{14}$$

and its density, which we denote as $(i^{(q)}(x), x \geq 0)$, by

$$i^{(q)}(x) = \sum_{n=0}^{\infty} \exp\left(-x/q^n\right) \frac{(-1)^n q^{\binom{n}{2}}}{(q; q)_\infty (q; q)_n}. \tag{15}$$

More precisely, (13) follows immediately from (12), whereas (14) can be derived from the factorization (5) of the exponential law. Finally, (15) is obtained by inversion of the Mellin transform; see [10] for details.

● **Brownian motion with drift:** Take $\xi_t = 2(B_t + bt)$, where $B$ is a standard Brownian motion and $b > 0$. Then $m = 2b$, $\psi(q) = 2q(q + b)$, and we get for every integer $k \geq 1$

$$\begin{aligned} m\frac{\psi(1)\cdots\psi(k-1)}{(k-1)!} &= 2^k b(b+1)\cdots(b+k-1) \\ &= 2^k \frac{\Gamma(k+b)}{\Gamma(b)}. \end{aligned}$$

The right-hand side can be identified as the $k$-th moment of $2\gamma_b$, where $\gamma_b$ is a gamma variable with index $b$. Hence Theorem 3 enables us to recover the identity in distribution

$$\int_0^\infty \exp\left\{-2(B_s + bs)\right\} ds \overset{\mathcal{L}}{=} \frac{1}{2\gamma_b}, \qquad (16)$$

which has been established by Dufresne [26] (see also Proposition 3 in Pollak-Siegmund [57], Example 3.3 on page 309 in Urbanik [65], Yor [68], ...). A further discussion of (16) in relation with DNA is made in [46] and [48].

• **Exponential Lévy measure:** Suppose $\xi$ has bounded variation with drift coefficient $+1$ and Lévy measure

$$\Pi(dx) = (a + b - 1)b\mathrm{e}^{bx}dx, \quad x < 0,$$

with $0 < a < 1 < a + b$. Then, $\psi(q) = q(q + 1 - a)/(b + q)$, $m = (1 - a)/b$, and we get for every integer $k \geq 1$

$$\begin{aligned}
m\frac{\psi(1)\cdots\psi(k-1)}{(k-1)!} &= \frac{(1-a)\cdots(k-a)}{b\cdots(k+b-1)} \\
&= \frac{\Gamma(k+1-a)\Gamma(b)}{\Gamma(1-a)\Gamma(k+b)}.
\end{aligned}$$

On the right-hand side, we recognize the $k$-th moment of a beta variable with parameter $(1 - a, a + b - 1)$. We thus get the identity in law

$$I_\infty \overset{\mathcal{L}}{=} \frac{1}{\beta_{1-a,a+b-1}},$$

which has been discovered by Gjessing and Paulsen [32].

### *4.2. An integro-differential equation*

In this section, we present an equation satisfied by the density of the law of the exponential functional. We shall assume here that the Lévy process is given in the form

$$\xi_t = ct + \sigma B_t + \tau_t^+ - \tau_t^-,$$

where $(B_t, t \geq 0)$ is a standard Brownian motion, $c$ a real number, $\tau^\pm$ two subordinators with no drift, such that $B$, $\tau^+$ and $\tau^-$ are independent. Using the same notation as in Theorem 1(vi), we shall further suppose that $\int_0^\infty \mathrm{e}^x \bar\Pi^-(x)dx < \infty$ and that

$$\int_0^\infty \bar\Pi^\pm(x)dx < \infty \quad, \quad c + \int_0^\infty \bar\Pi^+(x)dx - \int_0^\infty \bar\Pi^-(x)dx < \infty,$$

(which is equivalent to assuming that $\xi_1$ is integrable and has a positive mean, cf. the remark after Theorem 1). Carmona *et al.* [21] have proved the following result.

**Theorem 4** *Let the above assumptions prevail. Then the law of the terminal value $I_\infty$ of the exponential functional is absolutely continuous with respect to Lebesgue measure and has a $\mathcal{C}^\infty$ density $k$ which solves the following equation:*

$$-\frac{\sigma^2}{2}\frac{d}{dx}(x^2k(x)) + \left(\left(\frac{\sigma^2}{2}-c\right)x+1\right)k(x)$$

$$= \int_x^\infty \bar{\Pi}^+\left(\log(u/x)\right)k(u)du - \int_0^x \bar{\Pi}^-\left(\log(x/u)\right)k(u)du\,.$$

We refer to Carmona *et al.* [21] for a proof and for some applications of this equation.

### 4.3. Random affine equation

The strong Markov property of Lévy processes has an interesting consequence for the exponential functional. Indeed, consider some finite stopping time $T$ in the natural filtration $(\mathcal{F}_t)_{t\geq 0}$ of $\xi$, so that the process $\xi'_t := \xi_{T+t} - \xi_T$, $t \geq 0$, is independent of $\mathcal{F}_T$ and has the same law as $\xi$. Now writing

$$\int_0^\infty \exp(-\xi_s)ds \;=\; \exp(-\xi_T)\int_0^\infty \exp(-\xi'_r)dr + \int_0^T \exp(-\xi_s)ds$$

yields a so-called **random affine equation**

$$I_\infty \;\overset{\mathcal{L}}{=}\; AI_\infty + B \tag{17}$$

where, on the right-hand side, the variables $I_\infty$ and $(A,B)$ are independent. Specifically, the law of the latter pair is given by :

$$(A,B) \;\overset{\mathcal{L}}{=}\; \left(\exp(-\xi_T),\int_0^T \exp(-\xi_s)ds\right)\,.$$

The random affine equation has been investigated in depth by Kesten [45]; see also Vervaat [66] and Goldie [35].

An interesting application of the results of Kesten and Goldie in the special case when we take $T \equiv 1$ is the following estimation of the tail distribution of the exponential functional under a condition of Cramer's type.

**Corollary 5** *Suppose there exists $\theta > 0$ such that $\mathbb{E}(\exp(-\theta\xi_t)) = 1$ for some (and then all) $t > 0$. Assume further that $\mathbb{E}(\xi_t\exp(-\theta\xi_t)) < \infty$ and that $\xi$ is not arithmetic (i.e. there is no real number $r > 0$ such that $\mathbb{P}(\xi_t \in r\mathbb{Z}) = 1$ for all $t \geq 0$). Then there exists some constant $c > 0$ such that*

$$\mathbb{P}(I_\infty > t) \;\sim\; ct^{-\theta}\,, \qquad t \to \infty\,.$$

This interesting observation has been made by Méjane [53] and Rivero [58]. It constrasts with the situation when $\xi$ is a subordinator (then, obviously, the condition of Corollary 5 never holds), as we know that then $I_\infty$ possesses some finite exponential moments. Finally, we refer to the recent work [52] for extensions of Corollary 5 to situations where Cramer's condition is not fulfilled.

## 5. Self-similar Markov processes

Motivated by limit theorems for Markov processes, Lamperti [47] considered families of probability measures $(\mathbf{P}_x, x > 0)$ on Skorohod's space $\mathbb{D}$ of càdlàg paths $\omega : \mathbb{R}_+ \to \mathbb{R}_+$, under which the coordinate process $X.(\omega) = \omega(\cdot)$ becomes Markovian and fulfills the scaling property. That is for some real number $\alpha$, called the index of self-similarity, there is the identity

$$\text{the law of } (k X_{k^{-\alpha} t}, t \geq 0) \text{ under } \mathbf{P}_x \text{ is } \mathbf{P}_{kx} \tag{18}$$

where $k > 0$ and $x > 0$ are arbitrary. These laws are called self-similar (or semi-stable) Markov distributions.

Lamperti has shown that self-similar Markov processes can be connected to Lévy processes by a bijective correspondence that we shall present in the next section. Lamperti's correspondence also points at the key role of the exponential functional in this framework, and we shall see in particular in Section 5.3 that the law of $I_\infty$ is fundamental to understand the asymptotic behavior of self-similar Markov processes.

### 5.1. Lamperti's transformation

Recall that $\xi = (\xi_t, t \geq 0)$ denotes a real-valued Lévy process which starts from 0, and that the notation $\mathbb{P}, \mathbb{E}$ is used for probabilities and expectations related to $\xi$. In the sequel, it will be convenient to use the notation

$$I_t^{(b)} := \int_0^t \exp(-b \xi_s) ds \,,$$

for every $b \in \mathbb{R}$ (in other words, $I^{(b)}$ is the exponential functional associated to the Lévy process $b\xi$).

For each $\alpha \in \mathbb{R}$, we define a time-change $\tau^{(\alpha)} = (\tau_t^{(\alpha)}, t \geq 0)$ as the inverse of the exponential functional $I^{(-\alpha)}$, i.e.

$$\tau_t^{(\alpha)} := \inf\left\{ u \geq 0 : \int_0^u \exp(\alpha \xi_s) ds \geq t \right\} \,.$$

More precisely, the time-change is always finite when the Lévy process oscillates, or when $\xi$ drifts to $+\infty$ and $\alpha \geq 0$, or when $\xi$ drifts to $-\infty$ and $\alpha \leq 0$. In these cases, there is the identity

$$I_{\tau_t^{(\alpha)}}^{(-\alpha)} = \int_0^{\tau_t^{(\alpha)}} \exp(\alpha \xi_s) ds = t \,. \tag{19}$$

When $\xi$ drifts to $+\infty$ and $\alpha < 0$, or when $\xi$ drifts to $-\infty$ and $\alpha > 0$, (19) holds if and only if $t < I_\infty^{(-\alpha)}$, and $\tau_t^\alpha = \infty$ otherwise.

Next, for an arbitrary $x > 0$, denote by $\mathbf{P}_x$ the distribution on Skorokhod's space $\mathbb{D}$ of càdlàg paths of the process

$$X_t := x \exp\left\{\xi_{\tau_{tx^{-\alpha}}^{(\alpha)}}\right\}, \qquad t \geq 0.$$

It is this transformation from Lévy processes which we call Lamperti's transformation. In the case when the Lévy process drifts to $+\infty$ and $\alpha < 0$, we agree that the quantity above equals $\infty$ for all $t \geq x^\alpha I_\infty^{(-\alpha)}$. Similarly, we agree that $X_t = 0$ when the Lévy process drifts to $-\infty$, $\alpha < 0$ and $t \geq x^\alpha I_\infty^{(-\alpha)}$. In other words, $0$ and $\infty$ are the two absorbing boundaries; observe that these boundaries can only be reached continuously.

Then, writing $(\mathcal{G}_t)_{t \geq 0}$ for the natural filtration of the coordinate process $X$ on $\mathbb{D}$, we have that for every $s, t \geq 0$ and measurable function $f : \mathbb{R}_+ \to \mathbb{R}_+$

$$\mathbf{E}_x\left(f(X_{t+s}) \mid \mathcal{G}_t\right) = \mathbf{E}_y(f(X_s)) \qquad \text{for } y = X_t,$$

and the scaling property (18) holds. Conversely, any self-similar Markov process which either never reaches the boundary states $0$ and $\infty$, or reaches them continuously (in other words, there is no killing inside $]0, \infty[$) can be constructed by Lamperti's transformation from a Lévy process.

We further point that for every $\beta \neq 0$, the process $X$ raised to the power $\beta$, $X^\beta$, is self-similar with index $\alpha/\beta$; more precisely it can be constructed from Lamperti's transformation from the Lévy process $\beta\xi$. This enables to focus on the case when the index of self-similarity is $1$.

That exponential functionals play an important role in the study of self-similar Markov processes is plain in cases when absorption can occur, since the distribution of the absorption time under $\mathbf{P}_x$ is given by the identity

$$\zeta \overset{\mathcal{L}}{=} x^\alpha I_\infty^{(-\alpha)}.$$

We shall also see in the next sections that, at least in some special cases, there are explicit formulas for the moments of $X_t$ which bear similarities with those of the exponential functional, and that the terminal value of the exponential functional $I^{(\alpha)}$ is crucial to determine the asymptotic behavior of the self-similar Markov process.

### 5.2. Moments at a fixed time

Throughout this section, we shall assume that $\alpha$, the index of self-similarity of $X$, is positive.

First, we consider the case when the Lévy process $\xi$ is a subordinator; recall that $\Phi$ denotes its Laplace exponent, cf. formula (3). One gets

$$\mathbf{E}_x\left(X_t^{-p}\right) = x^{-p} \sum_{n=0}^{\infty} \frac{(-x^{-p}t)^n}{n!} \prod_{k=0}^{n-1} \Phi(p + k\alpha), \qquad (20)$$

where in the right-hand side, we agree that the product equals 1 for $n = 0$. See Equations (2) and (3) in [13] where this formula is proven for $x = \alpha = 1$ (the general case follows easily). It may be interesting to check that (20) agrees with Kolmogorov's backward equation. Specifically, if we write $\mathcal{G}$ for the infinitesimal generator of $X$, and $\varphi_p : x \to x^p$ for the power function with exponent $p \in \mathbb{R}$, then standard transformations for Markov processes (time substitution and image by a bijection) entail that

$$\mathcal{G}\varphi_{-p} = -\Phi(p)\varphi_{-p-\alpha}, \qquad p \geq 0.$$

It is then immediate to see that if we write $f_t(x)$ for the quantity in the right-hand side of (20), then $\partial_t f_t(x) = \mathcal{G}f_t(x)$ as required.

We also point out that in the special case $x = \alpha = p = 1$, formula (20) shows that the function $t \to \mathbf{E}_1\left(1/X_t\right)$ coincides with the Laplace transform of the probability measure $\rho_0$ which is related to the exponential functional in Theorem 2(ii).

Now we turn our attention the case when the Lévy process has finite exponential moments of any positive order, and which does not drift to $-\infty$. Specifically, we assume that (8) holds and that $m = \mathbb{E}(\xi_1) \geq 0$. Then, following [15], one can establish using Kolmogorov's backward equation that the entire moments of $X_t^\alpha$ are given for $k = 1, \ldots$ by

$$\mathbf{E}_x\left(X_t^{\alpha k}\right) = x^{\alpha k} + \sum_{\ell=1}^{k} \frac{\psi(\alpha k) \cdots \psi(\alpha(k - \ell + 1))}{\ell!} x^{\alpha(k-\ell)} t^\ell. \qquad (21)$$

Furthermore, if $\xi$ has no positive jumps, then (21) determines the semigroup of the self-similar Markov process.

### 5.3. Asymptotic behavior for positive indices of self-similarity

Throughout this section, we shall suppose that the index of self-similarity is positive : $\alpha > 0$, and that the Lévy process $\xi$ drifts to $+\infty$. In this situation, recall that the self-similar Markov process tends to $\infty$ without reaching $\infty$, and it is therefore interesting to obtain sharper information about the asymptotic behavior of $X_t$ as $t \to \infty$. The following result was proved in [14] (see also Bertoin and Caballero [11] in the special case when the Lévy process $\xi$ is a subordinator).

**Theorem 6** *Let $X$ denote a self-similar Markov process with index $\alpha > 0$ associated by Lamperti's transformation to the Lévy process $\xi$. Suppose that $\xi$ is not arithmetic (i.e. does not live on a discrete subgroup $r\mathbb{Z}$ for some $r > 0$), and has a positive finite mean, i.e. $\xi_1 \in L^1(\mathbb{P})$ and $\mathbb{E}(\xi_1) := m > 0$. Then for every $x > 0$, under $\mathbf{P}_x$, $t^{-1/\alpha}X_t$ converges in distribution as $t \to \infty$ to some probability measure $\mu_\alpha$ on $\mathbb{R}_+$ specified by*

$$\int_{\mathbb{R}_+} f(y)\mu_\alpha(dy) = \frac{1}{\alpha m}\mathbb{E}\left(\frac{1}{I_\infty^{(\alpha)}}f\left(\frac{1}{I_\infty^{(\alpha)}}\right)\right),$$

*where $f : \mathbb{R}_+ \to \mathbb{R}_+$ denotes a generic measurable function.*

**Example.** Suppose $\alpha = 1$ and $\xi_t = 2(B_t + bt)$ for some constant $b > 0$, i.e. $\xi$ is twice a standard Brownian motion with drift $b$. Then it is easy to check, e.g. by stochastic calculus, that $X$ solves the stochastic differential equation

$$dX_t = 2\sqrt{X_t}d\beta_t + 2(b+1)dt$$

where $(\beta_t, t \geq 0)$ is a standard Brownian motion. Thus $X$ is the square of a Bessel process with dimension $2(b+1)$. On the other hand, we know from Dufresne's identity (16) that $I_\infty^{(\alpha)} = I_\infty \overset{\mathcal{L}}{=} 1/2\gamma_b$ where $\gamma_b$ denotes a gamma variable with parameter $b$. Hence Theorem 6 yields that the limiting distribution $\mu_1$ is given by

$$\int_{\mathbb{R}_+} f(y)\mu_1(dy) = b^{-1}\mathbb{E}\left(\gamma_b f(2\gamma_b)\right) = \mathbb{E}(f(2\gamma_{b+1})),$$

i.e. $t^{-1}X_t$ converges in law as $t \to \infty$ towards twice a gamma variable with parameter $b + 1$.

Theorem 6 is the key to understand the entrance distribution of the self-similar Markov process from $0+$ (i.e. the limit of the laws $\mathbf{P}_x$ as $x \to 0+$); see [14] for details and also Caballero and Chaumont [19] for a recent development.

## 6. Some applications

Besides its role in the study of self-similar Markov processes, exponential functionals arise in a number of domains in probability theory and related topics, among which (the following list is by no mean exhaustive):

### 6.1. Brownian motion on hyperbolic spaces

The random variable

$$\int_0^\infty \exp(B_s^{(\nu)})dW_s^{(\mu)}, \tag{22}$$

where $W_s^{(\mu)} = W_s - \mu s$ and $B_s^{(\nu)} = B_s - \nu s$ are two independent Brownian motions with respective drifts $-\mu \in \mathbb{R}$ and $-\nu < 0$, appear both in risk theory (see e.g. Paulsen [55], Example 3.1) as well as in connection with invariant diffusions on the hyperbolic half-plane $\mathbb{H}$ (see Bougerol [16], Ikeda and Matsumoto [40], Baldi *et al.* [1], ...). Hence, it is no surprise that the distribution of (22) has been much studied; it has a density given by

$$f(x) = c_{\mu,\nu}\frac{\exp(-2\mu\arctan x)}{(1+x^2)^{\nu+1/2}}$$

which belongs to the type IV family of Pearson distributions.

Concerning diffusions on $\mathbb{H}$, consider the infinitesimal generator

$$L^{(\mu,\nu)} = \frac{y^2}{2}\Delta - \mu\nu\frac{\partial}{\partial x} - (\nu - 1/2)y\frac{\partial}{\partial y}$$

($(x,y)$ denote the Cartesian coordinates on $\mathbb{H}$). The differential operator $L^{(\mu,\nu)}$ is invariant under the orientation-preserving isometries of $\mathbb{H}$ that fix the point at infinity, that is under the real affine transformations $z \to az + b$ with $a > 0$ and $b \in \mathbb{R}$. The diffusion process associated to $L^{(\mu,\nu)}$ corresponds to the stochastic differential equation

$$\begin{cases} dX_t = Y_t dW_t - \mu Y_t dt \\ dY_t = Y_t dB_t - (\nu - 1/2)Y_t dt \,, \end{cases} \tag{23}$$

where, as before, $W$ and $B$ are independent one-dimensional Brownian motions. The solution of (23) with starting point $\mathtt{i}y = (0, y)$ is :

$$\begin{cases} X_t = y\int_0^t \exp(B_s^{(\nu)})dW_s^{(\mu)} \\ Y_t = y\exp(B_t^{(\nu)}) \,. \end{cases}$$

In close connection with this diffusion, one may derive a way to obtain the density $f(x)$ exhibited above. This may follow from the following three simple observations:

(a) For a fixed real number $x$, consider the two processes

$$X_t^{(\mu,\nu)} = \exp(B_t^{(\nu)})x + \int_0^t \exp(B_s^{(\nu)})dW_s^{(\mu)}\,,$$

and

$$\tilde{X}_t^{(\mu,\nu)} = \exp(B_t^{(\nu)})\left(x + \int_0^t \exp(-B_s^{(\nu)})dW_s^{(\mu)}\right)\,.$$

By time-reversal, these two processes have the same one-dimensional distributions, although they do not have the same law.

(b) The process $(\tilde{X}_t^{(\mu,\nu)}, t \geq 0)$ is a diffusion with generator

$$M = \frac{1 + x^2}{2}\frac{d^2}{dx^2} - (\mu + (\nu - 1/2)x)\frac{d}{dx}\,.$$

(c) The distribution at time $t$ of this diffusion process converges to the invariant distribution with density $f(x)$. Indeed, from (a), $X_t^{(\mu,\nu)}$ has the same law as $\tilde{X}_t^{(\mu,\nu)}$, and the former converges in law as $t \to \infty$ to, say, $X_\infty^{(\mu,\nu)}$. The limit distribution is invariant, that is annihilated by $M^*$, the adjoint of $M$, and hence it coincides with $f$.

## 6.2. Asymptotic studies of Brownian diffusions in random environments

Exponential functionals occur very naturally in the study of some models of Random Walks in Random Environment, which we now describe informally. Associated with a stochastic process $V = (V(x), x \in \mathbb{R})$ such that $V(0) = 0$, a diffusion $X_V = (X_V(t), t \geq 0)$ in the random potential $V$ is, loosely speaking, a solution to the stochastic differential equation

$$dX_V(t) = d\beta_t - \frac{1}{2}V'(X_V(t))dt, \qquad X_V(0) = 0,$$

where $(\beta_t, t \geq 0)$ is a standard Brownian motion independent of $V$. More rigorously, the process $X_V$ should be considered as a diffusion whose conditional generator, given $V$, is:

$$\frac{1}{2}\exp(V(x))\frac{d}{dx}\left(e^{-V(x)}\frac{d}{dx}\right).$$

It is now clear, from Feller's construction of such diffusions, that the potential $V$ does not need to be differentiable.

Brox [18] and Kawazu and Tanaka [42, 43] (see also e.g. Hu *et al.* [39] and Comtet *et al.* [23]) have studied this random diffusion in the case when $V$ is a Brownian motion with negative drift : $V(x) = B_x - kx$, for $k \geq 0$. Thus, it is no surprise that the knowledge about the exponential functionals $\int_0^x dy \exp(B_y - ky)$, $x \in \mathbb{R}$, plays an essential role in this domain; see Shi [61] for a nice picture of this area.

## 6.3. Mathematical finance

In mathematical finance, the processes of reference were chosen originally to be the exponential of Brownian motions with drift (i.e. the famous Black & Scholes framework, which has been widely adopted between 1973 and 1990, say), and later generalized to become the exponential of Lévy processes ($\xi_t, t \geq 0$). In particular, the computation of the price of Asian options (see, e.g., papers in [70]) is equivalent to the knowledge of the law of $\int_0^t \exp(-\xi_s)ds$, for fixed $t$; an important reduction of the problem consists in replacing $t$ by an exponential time $\theta$, which is independent of $\xi$ (in other words, one kills the Lévy process with some constant rate). Furthermore, some identities in law due to Dufresne [27], for exponential functionals of Brownian motions with opposite drifts, involving both the case $t < \infty$ and $t = \infty$, led to some interesting extensions of Pitman's $(2M - X)$ theorem for exponential functionals (see, e.g. [49, 50, 51]).

## 6.4. Self-similar fragmentations

Fragmentations are meant to describe the evolution of an object with unit mass which falls down randomly into pieces as time passes. They arise in a variety of models; to give just a few examples, let us simply mention the studies of stellar fragments and meteoroids in astrophysics, fractures and earthquakes in geophysics, breaking of crystals in crystallography, degradation of large polymer chains in chemistry, fission of atoms in nuclear physics, fragmentation of a hard drive in computer science, crushing in the mining industry, ...

The state of the system at some given time consists in the sequence of the masses of the fragments. Suppose that its evolution is Markovian and obeys the following rule. There is a parameter $\alpha \in \mathbb{R}$, called the index of self-similarity, such that given that the system at time $t \geq 0$ consists in the ranked sequence of masses $m_1 \geq m_2 \geq \ldots \geq 0$, the system at time $t + r$ is obtained by dislocating every mass $m_i$ independently of the other fragments to obtain a family of sub-masses, say $(m_{i,j}, j \in \mathbb{N})$, where the sequence of the ratios $(m_{i,j}/m_i, j \in \mathbb{N})$ has the same distribution as the sequence resulting from a single unit mass fragmented up to time $m_i^\alpha r$. Such a random process will be referred to as a *self-similar fragmentation* with index $\alpha$.

For the sake of simplicity, we shall focus here on the following elementary model. A fragment with mass $m > 0$ remains unchanged during an exponential time with parameter $m^\alpha$. Then it splits and gives rise to a family of smaller masses, $ms_1, ms_2, \ldots$ where, $(s_1, s_2, \ldots)$ is a random sequence of nonnegative integers mass $m_i$ with sum $s_1 + s_2 + \ldots = 1$, which has a fixed distribution $\nu$. One calls $\nu$ the dislocation law of the fragmentation process.

In this setting, self-similar Markov processes arise as follows. One taggs a point at random in the object according to the mass measure and independently of the fragmentation process, and one considers at each time $t$ the mass $\chi_t$ of the fragment containing the tagged point. Then $(1/\chi_t, t \geq 0)$ is a self-similar Markov process with index $\alpha$, and the Lévy measure $\Pi$ of the Lévy process $\xi$ which is associated to $1/\chi$ by Lamperti's transformation (cf. Section 5.1), can be expressed explicitly in terms of the dislocation law $\nu$ as

$$\Pi(dx) = e^{-x} \sum_{j=1}^{\infty} \nu(-\log s_j \in dx), \qquad x \in ]0, \infty[.$$

Statistical properties of the tagged fragment can then be derived e.g. from results presented in Section **??**, and then be shifted to the fragmentation process. In particular, this approach yields important limit theorems for the empirical distribution of the fragments of self-similar fragmentations with a positive index of self-similarity; see [9, 12, 17, 30]. On the other hand, self-similar fragmentations with a negative index of self-similarity are dissipative, in the sense that the total mass of the fragments decreases as time passes. In this direction, the terminal value of the exponential functional gives the absorption time of the self-similar tagged process, i.e. the time at which the mass of the tagged fragment vanishes. Distributional properties of the exponential functional provide a crucial tool for the analysis of the loss of mass; see Haas [37, 38].

In a related area, we mention that the terminal value of the exponential functional of subordinators also arises in the study of regenerative compositions, see the recent work of Gnedin, Pitman and Yor [33, 34].